\documentclass{amsart}
\usepackage{amssymb}
\title{Proving Touchard's Theorem From Euler's Form }
\author{Eyob Delele Yirdaw}

\begin{document} 

\address{Addis Ababa University,
Addis Ababa, Ethiopia}
\email{edyirdaw@sisa.aau.edu.et}

\subjclass[2000]{Primary 11A25}

\date{March 18, 2008.}

\keywords{Odd perfect numbers, Touchard's theorem}
\begin{abstract}
This paper derives Touchard's theorem from Euler's form for odd perfect numbers. It also fine-tunes Euler's form.
\end{abstract}
\maketitle
A perfect number is a natural number $N$ which satisfies $\sigma(N)=2N$, where $\sigma(N)$ is defined as the sum of the divisors of $N$. Because $\sigma(ab)=\sigma (a)\sigma(b)$ when $gcd(a,b)=1$, $\sigma(N)$ is a multiplicative function. A handful of even perfect numbers are known and it is believed that infinite of them exist. It is also believed that odd perfect numbers do not exist. Euler showed that an odd perfect number $N$ must be of the form $N=p^\alpha Q^2$ where $p$ is a prime, $p\equiv\alpha\equiv 1\pmod{4} $ and $Q$ is an odd number such that $p\nmid Q$. Touchard \cite{A} showed that $N=12k+1$ or $N=36k+9$. Euler's form manages to imply only $N\equiv 1 \pmod{4}$. We have to refine it a little bit to get to Touchard's theorem. Since $p$ is prime, it can have only the form $12n+1$ or $12n+5$. Because $\alpha$ is always odd, $(p+1)\mid \sigma(p^\alpha)$. Taking $p=12n+5$, we have $(12n+6)\mid\sigma((12n+5)^\alpha)$. Because $3\mid(12n+6)$, we have $3\mid N$. Euler's form then becomes
\begin{equation}\label{eq-1}
N=(12n+1)^\alpha Q^2 \qquad or \qquad N=(12n+5)^\alpha (3Q)^2.
\end{equation}
Using (1) we can prove Touchard's theorem. When $3\mid N$, (1) becomes
\begin{eqnarray}
N=p^\alpha(3Q)^2=9.p^\alpha Q^2=9(4k+1)=36k+9. \nonumber
\end{eqnarray}
When $3\nmid N$, we use only $N=(12n+1)^\alpha Q^2$, $3\nmid Q$ (i.e., $Q \equiv 1\pmod{6}$ or $Q \equiv 5\pmod{6}$). Because $Q^2 \equiv 1\pmod{12}$, $N=(12n+1)^\alpha Q^2 \equiv 1\pmod{12}$, completing the proof.

We can refine (1) a little bit further to 
\begin{align}
N&=(12n+1)^{12\lambda+1}Q^2  \qquad or \nonumber \\
N&=(12n+1)^{12\lambda+9}Q^2  \qquad or \nonumber \\
N&=(12n+5)^{12\lambda+1}(3Q)^2 \qquad or \nonumber \\
N&=(12n+5)^{12\lambda+9}(3Q)^2 \qquad or  \\
N&=(12n+1)^{12\lambda+5}(3Q)^2  \quad for \quad n\equiv 0,1\pmod{7} \quad or \nonumber \\
N&=(12n+5)^{12\lambda+5}(3Q)^2  \quad for \quad n\equiv 2,3\pmod{7} \quad or \nonumber \\
N&=(12n+1)^{12\lambda+5}(21Q)^2  \quad for \quad n\equiv 2,3,5,6\pmod{7} \quad or \nonumber \\
N&=(12n+5)^{12\lambda+5}(21Q)^2  \quad for \quad n\equiv 0,1,4,5 \pmod{7}. \nonumber 
\end{align}

When $p=12n+1$ and $\alpha=12\lambda+5$,
\begin{equation}
\sigma((12n+1)^{12\lambda+5})=[(12n+1)+1][(12n+1)^4+(12n+1)^2+1]\sum_{\gamma=0}^{2\lambda}
(12n+1)^{6\gamma}. 
\end{equation}
Here $3\mid [(12n+1)^4+(12n+1)^2+1]$. Thus $3\mid N$ and (1) becomes $N=(12n+1)^{12\lambda+1} Q^2$ or $ N=(12n+1)^{12\lambda+9} Q^2$ or $N=(12n+1)^{12\lambda+5} (3Q)^2$ or $ N=(12n+5)^{4\lambda +1}(3Q)^2$. If $n$ in (3) is congruent to $2,3,5$ or $6$ mod(7), then $7\mid [(12n+1)^4+(12n+1)^2+1]$ and hence $7\mid N$. When $n\equiv 4\pmod{7}$ $p$ will not be prime. When $p=12n+5$ and $\alpha=12\lambda+5$, an argument similar to the above one gives $7\mid N$ for $n\equiv 0,1,4,5 \pmod{7}$.

Kuhnel's theorem \cite{C} that $105\nmid N$ requires $5\nmid N$ in the last two forms of (2).

\bibliographystyle{amsplain}

\end{document}